\documentclass[12pt]{amsart}
\usepackage{amscd}
\theoremstyle{definition}

\def \ph{\varphi}

\def \refeq#1{equation (\ref{#1})}
\def \ra{\rightarrow}

\def \hom{\mbox{\rm Hom}}

\def \tns{\otimes}

\def \k{\mbox{$\mathfrak K$}}

\def \C{\mbox{$\mathbb C$}}
\def \Z{\mbox{$\mathbb Z$}}

\def\br#1#2{\lbrack#1,#2\rbrack}

\def\zt{\mbox{$\Z_2$}}

\def\sh{\operatorname{Sh}}

\def\inv{^{-1}}
\def\d{d}
\def\td{\tilde\d}

\def\A{\mbox{$\mathcal A$}}
\def\B{\mbox{$\mathcal B$}}
\def\L{L}

\def\LA{\mbox{$\L_{\A}$}}
\def\LB{\mbox{$\L_{\B}$}}

\def\m{\mbox{$\mathfrak m$}}
\def\a{\mbox{$\mathfrak a$}}

\def\coder{\operatorname{Coder}}

\def\linf{\mbox{$L_\infty$}}
\def\and{\mbox{ \rm and }}

\def\s#1{(-1)^{#1}}
\DeclareMathOperator*{\invlim}{\overleftarrow{\rm lim}}

\def\htns{\hat\tns}
\def\SW{S(W)}

\def\htns{\hat\tns}

\def\two{e_2}
\def\one{e_1}

\def\inv{^{-1}}
\def\one{\mathbf{1}}
\def\two{\mathbf{2}}
\def\tre{\mathbf{3}}
\def\dt#1{d(#1)}
\author{Alice Fialowski}
\address{E\"otv\"os Lor\'and University\\
Department of Applied Analysis\\
H-1117 Budapest, P\'azm\'any P. s\'et\'any. 1/C, HUNGARY}
\email{fialowsk@cs.elte.hu}
\author{Michael Penkava}
\address{University of Wisconsin\\
Department of Mathematics\\
Eau Claire, WI 54702-4004} \email{penkavmr@uwec.edu}
\thanks{The research of the authors was partially supported by
grants MTA-OTKA-NSF 38453, OTKA T043641, T043034, and
 grants from the University of Wisconsin-Eau Claire}
\title[Deformations of 3-d Lie Algebras]{Versal Deformations of Three
Dimensional Lie algebras as \linf\ Algebras}

\begin{document}
\setlength{\multlinegap}{0pt}
%\nocite{ps2,pen1,ls,pen2,pen3,kon,
%fm,mar,mar2,ksv,sta1,sta2,sta3,getz,getz2,ge_ka1,ge_ka2,
%gers,lod,hoch,fi,fi2,ff2,ff3}
\begin{abstract}
We consider versal deformations of
$0|3$-dimensional \linf\ algebras,  which correspond precisely to
ordinary (non-graded) three dimensional Lie algebras. The
classification of such algebras over $\C$ is well known, although
we shall give a derivation of this classification using an
approach of treating them as \linf\ algebras. Because the
symmetric algebra of a three dimensional odd vector space contains
terms only of exterior degree less than or equal to three, the
construction of versal deformations can be carried out completely.
We give a
characterization of the moduli space of Lie algebras using
deformation theory as a guide to understanding the picture.
\end{abstract}

\date\today
\maketitle
%\table
%\input intro.tex
\section{Introduction}

The classification of low dimensional Lie algebras has been known
for a long time. For example, the classification of ordinary Lie
algebras of dimension 3, the subject of this paper, appears in
textbooks such as \cite{jac}. More recently, the moduli space of
three dimensional Lie algebras was studied in \cite{aga,tu}. The
problem of finding a versal deformation of a given object is a
basic question in deformation theory because such a deformation
induces all other deformations.  This problem turns out to be very
difficult. Versal deformation theory was first worked out for the
case of Lie algebras in  \cite{fi, ff2} and then extended to
$L_\infty$ algebras in \cite{fp1}.  We apply these general results
to construct versal deformations of three dimensional ordinary Lie
algebras, treating them as examples of \linf\ algebras. We use the
methods developed in \cite{fp1, fp2}.

\linf\ algebras are natural generalizations of Lie algebras and
superalgebras if one considers $\Z_2$-graded vector spaces. An ordinary
$3$-dimensional Lie algebra is the same thing as a \linf\
algebra structure on a $0|3$ (0 even and 3 odd) dimensional \zt-graded
vector space.  \linf\ algebras
were first described in \cite{ss} and have recently been the focus
of much attention \cite{lm,bfls,ls,sta3}. The advantage of
considering Lie algebras as \linf\ algebras is that the
deformation problem becomes simpler and we get a clearer insight
to the moduli space of the variety of Lie algebras in a given
dimension.

Even though \linf\ algebras appear in many contexts, there are
only a few examples known, even in low dimensions. In \cite{fp2},
the authors classified \linf\ algebras of dimension less than or
equal to 2. In related work \cite{fp4}, we studied \linf\ algebras
of dimension $1|2$, and in recent work with Derek Bodin
\cite{bfp1}, we studied \linf\ algebras of dimension $2|1$.  Our
goal in this paper is to study versal deformations of $0|3$
dimensional \linf\ algebras.  The examples we study here are much
simpler than the ones studied in \cite{fp4, bfp1}, so we are able
to give a complete treatment of versal deformations.

For simplicity we will suppose that the underlying vector space is
defined over $\C$. Note that in the classification of Lie algebra
structures, this assumption reduces the number of equivalence
classes of Lie algebra structures. For example,  over the reals,
there are two non-equivalent classes of simple Lie algebra
structures in dimension 3, while over \C, there is only one class.
However,  we shall not consider the more general question of
deformations over other fields.

In Section~2 we introduce $L_\infty$ algebras and give the definition
of a versal deformation --- indicating a construction we will use in
our computation.  Section~3 treats the classification of
codifferentials, giving another approach to the well-known
classification of  $3$-dimensional Lie algebras.  In Section~4 we
compute versal deformations of each of these Lie algebras and also,
using our deformation results, describe their moduli space.

\section{Basic Definitions}
\subsection{ $L_\infty$ Algebras}
If $W$ is a \zt-graded vector space, then $\SW$ denotes the symmetric
coalgebra of $W$. If we let $T(W)$ be the reduced tensor algebra
$T(W)=\bigoplus_{n=1}^\infty W^{\tns n}$, then the reduced symmetric algebra
$S(W)$ is the quotient of the tensor algebra by the graded ideal
generated by $u\tns v-\s{uv}v\tns u$ for elements $u$, $v\in W$.
The symmetric algebra has a natural coalgebra structure, given by
\begin{equation*}
\Delta(w_1\dots w_n)=
\sum_{k=1}^{n=1}\sum_{\sigma\in\sh(k,n-k)}\epsilon(\sigma)
w_{\sigma(1)}\dots w_{\sigma(k)}\tns w_{\sigma(k+1)}\dots w_{\sigma(n)},
\end{equation*}
where we denote the product in $S(W)$ by juxtaposition,
$\sh(k,n-k)$ is the set of \emph{unshuffles} of type $(k,n-k)$,
and $\epsilon(\sigma)$ is a sign determined by $\sigma$ (and $w_1\dots w_n$)
given by
$$
w_{\sigma(1)}\dots w_{\sigma(n)}=\epsilon(\sigma)w_1\dots w_n.
$$
A {\it coderivation} on $S(W)$ is a map $\delta:S(W)\ra S(W)$ satisfying
$$\Delta\circ\delta=(\delta\tns I+I\tns\delta)\circ\Delta.$$
Let us suppose that the even part of $W$ has basis
$e_1\dots e_m$, and the odd part has basis $f_1\dots f_n$, so that $W$
is an $m|n$ dimensional space. Then a basis of $\SW$ is given by all
vectors of the form $e_1^{k_1}\dots e_m^{k_m}f_1^{l_1}\dots f_n^{l_n}$,
where $k_i$ is any nonnegative integer, and $l_i\in\zt$. An \linf\
structure on $W$ is simply an odd codifferential on $\SW$, that is to
say, an odd coderivation whose square is zero. The space $\coder(W)$
can be naturally identified with $\hom(S(W),W)$, and the Lie
superalgebra structure on $\coder(W)$ determines a Lie bracket on
$\hom(S(W),W)$ as follows. Denote $L_m=\hom(S^m(W),W)$ so that
$L=\hom(S(W),W)$ is the direct product of the spaces $L_i$.
If $\alpha\in L_m$ and $\beta\in L_n$, then $\br\alpha\beta$ is the
element in $L_{m+n-1}$ determined by
\begin{multline}\label{braform}
\br\alpha\beta
(w_1\dots w_{m+n-1})=\\
\sum_{\sigma\in\sh(n,m-1)}\epsilon(\sigma)\alpha(\beta(w_{\sigma(1)}\dots
w_{\sigma(n)})w_{\sigma(n+1)}\dots w_{\sigma(m+n-1)})\\
-\s{\alpha\beta}\
\sum_{\sigma\in\sh(m,n-1)}\epsilon(\sigma)
\beta(\alpha(w_{\sigma(1)}\dots
w_{\sigma(m)})w_{\sigma(m+1)}\dots w_{\sigma(m+n-1)}).
\end{multline}
Another way to express this bracket is in the form
\begin{equation*}
\br\alpha\beta=\alpha\tilde\beta-\s{\alpha\beta}\beta\tilde\alpha,
\end{equation*}
 where for $\ph\in\hom(S^k(W),W)$, $\tilde\ph$ is the
associated coderivation, given by
\begin{equation*}
\tilde\ph(w_1\dots w_n)=
\sum_{\sigma\in\sh(k,n-k)}\epsilon(\sigma)
\ph(w_{\sigma(1)}\dots w_{\sigma(k)})w_{\sigma(k+1)}\dots w_{\sigma(n)}.
\end{equation*}

If $W$ is completely odd, and $d\in\L_2$, then $d$ determines an ordinary
Lie algebra on $W$, or rather on its parity reversion.  This is the case
we consider in the present paper.  The
symmetric algebra on $W$ looks like the exterior algebra on $W$ if
we forget the grading.  If we define $[a,b]=d(ab)$ for $a,b\in W$, then
the bracket is antisymmetric because $ba=-ab$, and moreover
\begin{align*}
0=\br dd(abc)=&\
\frac12\sum_{\sigma\in\sh(2,1)}\epsilon(\sigma)
d(d(\sigma(a)\sigma(b))\sigma(c))\\=&
d((d(ab)c)+d(d(bc)a)-d(d(ac)b)\\
=&[[a,b],c]+[[b,c],a]-[[a,c],b],
\end{align*}
which is the Jacobi identity.  When $d\in\L_2$ and $W$ has a true grading,
then the same principle holds, except that one has to take into
account a sign arising from the map $S^2(W)\ra\bigwedge^2(V)$, where
$V$ is the parity reversion of $W$.  Thus \zt-graded Lie algebras
are also examples of \linf\ algebras.  In addition, differential
graded Lie algebras, also called superalgebras, are examples of \linf\ algebras.  In all these
cases, the method of construction  of miniversal deformations we
describe here applies.  One simply considers only terms that
come from $\L_2$ in the codifferentials.

Suppose that $\tilde g:S(W)\ra S(W')$ is a coalgebra morphism,
that is a map satisfying
$$
\Delta'\circ\tilde g=(\tilde g\tns \tilde g)\circ\Delta.
$$
If $d$ and $d'$ are \linf\ algebra structures on $W$ and $W'$,
resp., then $\tilde g$ is a homomorphism between these structures if
$\tilde g\circ d=\d'\circ\tilde g$. Two \linf\ structures $d$ and $d'$ on $W$
are equivalent, and we write $d'\sim d$ when there is a coalgebra automorphism $\tilde g$ of $S(W)$
such that $d'={\tilde g}^*(d)={\tilde g}\inv\circ d\circ {\tilde g}$.   Furthermore,
if $d=d'$, then $\tilde g$ is said to be an automorphism of the \linf\
algebra.

\subsection{ Versal Deformations}

An augmented local ring $\A$ with maximal ideal $\m$ will be called an
\emph{infinitesimal base} if  $\m^2=0$, and a  \emph{formal base} if $\A=\invlim_n \A/\m^n$.
A deformation of an \linf\ algebra structure $d$ on $W$
with  base given by a local ring $\A$ with augmentation $\epsilon:\A\ra\k$, where
$\k$ is the field over which $W$ is defined, is an $\A$-\linf\ structure $\td$ on $W\htns\A$
such that the morphism of $\A$-\linf\ algebras
$\epsilon_*=1\tns\epsilon:\LA=L\tns \A\ra L\tns\k=L$ satisfies
$\epsilon_*(\td)=d$.  (Here $W\htns\A$ is an appropriate completion of $W\tns\A$.)
The deformation is called infinitesimal (formal) if $\A$ is an infinitesimal (formal) base.

In general, the cohomology $H(D)$ of $d$ given by the operator
$D:L\ra L$ with $D(\ph)=\br\ph d$ may not be finite dimensional.
However, $L$ has a natural
 filtration $L^n=\prod_{i=n}^\infty L_i$,
which induces a filtration $H^n$ on the cohomology, because
 $D$
respects the filtration. Then $H(D)$ is of finite type if $H^n/H^{n+1}$
is finite dimensional.
 Since this is always true when $W$ is finite
dimensional, the examples we study here will always
 be of finite
type.  A set $\delta_i$ will be called a \emph{basis of the cohomology},
if any element $\delta$
 of the cohomology can be expressed uniquely
as a formal sum $\delta=\delta_i a^i$. If we identify
 $H(D)$ with a
subspace of the space of cocycles $Z(D)$, and we choose a \emph{basis}
$\beta_i$ of the
 coboundary space $B(D)$, then any element
$\zeta\in\Z(D)$ can be expressed uniquely as a sum
$\zeta=\delta_i a^i +\beta_i b^i$.

For each $\delta_i$, let $u^i$ be a parameter of opposite parity.  Then the infinitesimal
deformation $d^1=d+\delta_i u^i$, with base $\A=\k[u^i]/(u^iu^j)$ is universal in the sense that
if $\td$ is any infinitesimal deformation with base $\B$, then there is a unique homomorphism
$f:\A\ra\B$, such that the morphism $f_*=1\tns f:\LA\ra\LB$ satisfies $f_*(\td)\sim d$.

For formal deformations, there is no universal object in the sense above. A
\emph{versal deformation}
is a deformation $d^\infty$ with formal base $\A$ such that if $\td$ is any formal deformation
with base $\B$, then there is some morphism $f:\A\ra\B$ such that $f_*(d^\infty)\sim\td$. If $f$
is unique whenever $\B$ is infinitesimal, then the versal deformation is called \emph{miniversal}.
In \cite{fp1}, we constructed a miniversal deformation for \linf\ algebras with finite type
cohomology.

The method of construction is as follows.  Define a coboundary
operator $D$ by $D(\ph)=[\ph,d]$. First, one constructs the universal
infinitesimal deformation $d^1=d+\delta_i u^i$, where $\delta_i$
is a graded basis of the cohomology $H(D)$ of $d$, or more
correctly, a basis of a subspace of the cocycles which projects
isomorphically to a basis in cohomology, and $u^i$ is a parameter
whose parity is opposite to $\delta_i$.  The infinitesimal
assumption that the products of parameters are equal to zero gives
the property that $[d^1,d^1]=0$. Actually, we can express
\begin{equation*}
[d^1,d^1]=\s{\delta_j(\delta_i+1)}[\delta_i,\delta_j]u^iu^j=\delta_k
a^k_{ij}u^iu^j +\beta_k b^k_{ij}u^iu^j,
\end{equation*}
where $\beta_i$ is a basis of the coboundaries, because the
bracket of $d^1$ with itself is a cocycle. Note that the right
hand side is of degree 2 in the parameters, so it is zero up to
order 1 in the parameters.

If we suppose that $D(\gamma_i)=-\frac12\beta_i$, then by
replacing $d^1$ with $$d^2=d^1+\gamma_kb^k_{ij}u^iu^j,
$$ one obtains
\begin{equation*}
[d^2,d^2]=\delta_k a^k_{ij}u^iu^j+ 2[\delta_l
u^l,\gamma_kb^k_{ij}u^iu^j]+[\gamma_kb^k_{ij}u^iu^j,\gamma_lb^l_{ij}u^iu^j]
\end{equation*}
Thus we are able to get rid of terms of degree 2 in the
coboundary terms $\beta_i$, but those which involve the cohomology
terms $\delta_i$ can not be eliminated.  This gives rise to a set
of second order relations on the parameters. One continues this
process, taking the bracket of the \emph{$n$-th order deformation}
$d^n$, adding some higher order terms to cancel coboundaries,
obtaining higher order relations, which extend the second order
relations.

Either the process continues indefinitely, in which
case the miniversal deformation is expressed as a formal power
series in the parameters, or after a finite number of steps, the
right hand side of the bracket is zero after applying the $n$-th
order relations. In this case, the miniversal deformation is
simply the $n$-th order deformation. In any case, we obtain a set
of relations $R_i$ on the parameters, one for each $\delta_i$, and
the algebra $A=\C[[u^i]]/(R_i)$ is called the base of the
miniversal deformation. Examples of the construction of miniversal
deformations can be found in \cite{ff3,ff2,fp2}.

\section{Classification of Lie Algebra Structures of Dimension~3}
Suppose that $W=\langle f_1, f_2, f_3 \rangle$.  Then $S(W)$
decomposes into three pieces.
\begin{equation*}
\begin{matrix}
S^1(W)&=\langle f_1,f_2,f_3 \rangle,&\qquad \dim(S^1(W))=0|3\\
S^2(W)&=\langle f_1f_2,f_1f_3,f_2f_3\rangle,&\qquad \dim(S^2(W))=3|0\\
S^3(W)&=\langle f_1f_2f_3\rangle,&\qquad \dim(S^3(W))=0|1.
\end{matrix}
\end{equation*}
Let $L=\hom(S(W),W)$ and $L_n=\hom(S^n(W),W)$.  Then
\begin{equation*}
\begin{matrix}
L_1(W)&=&\{\ph^I_j|\- I\in\{100,010,001\},j=1\dots 3\},&\quad \dim(L_1)=9|0\\
L_2(W)&=&\{\ph^I_j|I\in\{110,101,011\},j=1\dots 3\},&\quad \dim(L_2)=0|9\\
L_3(W)&=&\{\ph^{111}_j|j=1\dots 3\},&\quad \dim(L_3)=3|0\\
\end{matrix}
\end{equation*}
It follows that the only candidate for an odd codifferential is of
the form
\begin{align}
\notag d=&\ph^{110}_1a_1+\ph^{110}_2a_2+\ph^{110}_3a_3\\+
  &\ph^{101}_1a_4+\ph^{101}_2a_5+\ph^{101}_3a_6\\\notag+
  &\ph^{011}_1a_7+\ph^{011}_2a_8+\ph^{011}_3a_9
\end{align}
Being a quadratic codifferential,  we see that $d$ gives an \linf\
structure precisely when it determines a Lie algebra structure.
It is natural to consider the derived subalgebra $W'=d(S^2(W))$.
Let
$$
A=\begin{pmatrix}a_1&a_2&a_3\\a_4&a_5&a_6\\a_7&a_8&a_9\end{pmatrix}.
$$
It is easy to see that the rank of $A$ is precisely equal to the
dimension of the derived subalgebra. In particular,  when
$\det(A)=0$, the derived subalgebra has dimension less than three.

The codifferential condition $[d,d]=0$ is equivalent to the system
of three quadratic equations
\begin{align}
\notag
a_9a_2-a_3a_8+a_6a_1-a_3a_4&=0\\
-a_5a_9+a_8a_6+a_5a_1-a_2a_4&=0\\
\notag -a_4a_9-a_1a_8+a_7a_6+a_7a_2&=0.
\end{align}
Letting $x=a_2+a_6$, $y=a_9-a_1$, and $z=-(a_4+a_8)$, we can rewrite
the above equations as homogeneous linear equations for $x$, $y$, $z$:
$$
\begin{pmatrix}a_1&a_2&a_3\\a_4&a_5&a_6\\a_7&a_8&a_9\end{pmatrix}
\begin{pmatrix}x\\y\\z\end{pmatrix}=0.
$$
If $\det\,A \ne 0$ then the only possibility is $x=0$, $y=0$, and
$z=0$, that is, $a_6=-a_2$, $a_9=a_1$, and $a_8=-a_4$.  Thus
\begin{equation}\label{3ddet}
A=\begin{pmatrix}a_1&a_2&a_3\\a_4&a_5&-a_2\\a_7&-a_4&a_1\end{pmatrix},
\end{equation} whose determinant
$\det\,A =a_5a_1^2-2a_1a_2a_4-a_3a_4^2-a_7a_2^2-a_7a_3a_5,$ does
not vanish in general.  Consequently we have only one pattern to
consider for when the derived subalgebra has dimension three,
\begin{align}\label{3dderived}
\notag d=&\ph^{110}_1a_1+\ph^{110}_2a_2+\ph^{110}_3a_3\\+
  &\ph^{101}_1a_4+\ph^{101}_2a_5-\ph^{101}_3a_2\\\notag+
  &\ph^{011}_1a_7-\ph^{011}_2a_4+\ph^{011}_3a_1
\end{align}

When the derived subalgebra has dimension one, we can choose a
basis such that the codifferential $d$ has the simple form
\begin{equation}\label{1dderived}
d=\ph^{110}_1a_1+\ph^{101}_1a_4+\ph^{011}_1a_7.
\end{equation}
Moreover, it is easy to check that any coderivation of this form
is a codifferential.

Finally, suppose that the derived subalgebra has dimension 2.
Then we can express $d$ in the form
\begin{equation*}
d=\ph^{110}_1a_1+\ph^{110}_2a_2+\ph^{101}_1a_4+\ph^{101}_2a_5+\ph^{011}_1a_7+\ph^{011}_2a_8.
\end{equation*}
However, in this case, it is easy to check that the condition
$d^2=0$ forces $a_1=0$ and $a_2=0$, or the rank of the associated
matrix $A$ must be less than 2. Thus, when the derived subalgebra
has dimension 2, in an appropriate basis, we can express
\begin{equation}\label{2dderived}
d=\ph^{101}_1a_4+\ph^{101}_2a_5+\ph^{011}_1a_7+\ph^{011}_2a_8.
\end{equation}

Each of these three cases can be reduced to a much simpler form, up
to equivalence. Let us begin with the one dimensional case first.
Let us suppose that $g$ is a linear automorphism of the symmetric
coalgebra of $W$, in other words, $g$ is generated by an
invertible linear operator on $W$. For convenience,  let us denote
$f_1=\one$, etc. Suppose we express
\begin{align}
\notag
g(\one)&=\one l+\two m+\tre n\\
g(\two)&=\one r+\two s+\tre t\\
\notag g(\tre)&=\one x+\two y+\tre z
\end{align}
Let $G=\begin{pmatrix}l&m&n\\r&s&t\\x&y&z\end{pmatrix}$ be the
matrix of the linear transformation, so we require $\det(G)\ne0$.
Let us note that on $S^2(W)$, $g$ is given by
\begin{align*}
g(\one\two)&=\one\two(ls-mr)+\one\tre(lt-nr)+\two\tre(mt-ns)\\
g(\one\tre)&=\one\two(ly-mx)+\one\tre(lz-nx)+\two\tre(mz-ny)\\
g(\two\tre)&=\one\two(ry-sx)+\one\tre(rz-tx)+\two\tre(sz-ty)
\end{align*}

Now two codifferentials $d'$ and $d$ are equivalent precisely when
there is some automorphism $g$ such that $d'=gdg\inv$, in other
words, $d'g=gd$.

Let $d=\ph^{110}_1a_1+\ph^{101}_1a_4+\ph^{011}_1a_7$, and
$d'=\ph^{011}_1$.

Examining the conditions for these two codifferentials to be
equivalent, we obtain
\begin{align*}
d'g(\one\two)&=\one(mt-ns)\\
gd(\one\two)&=\one la_1+\two ma_1+\tre na_1\\
d'g(\one\tre)&=\one(mz-ny)\\
gd(\one\tre)&=\one la_4+\two ma_4+\tre na_4\\
d'g(\two\tre)&=\one(sz-ty)\\
gd(\two\tre)&=\one la_7+\two ma_7+\tre na_7\\
\end{align*}
Note that if all three coefficients $a_i$ vanished, this would
force $\det(G)=0$, which is not permitted. Thus we must have
$m=n=0$. But then it follows that $a_1=a_3=0$.  This means the only
members of the family of one dimensional codifferentials
equivalent to $d'$ are the nonzero multiples of $d'$.

Let us now consider the case when $d'=\ph^{110}_1$.  By a similar
analysis to the above, it is easily shown that if either $a_1$ or
$a_3$ does not vanish, then $d'$ and $d$ are equivalent.  Thus all
one dimensional solutions are equivalent to one of two
codifferentials $d=\ph^{110}_1$ or $d=\ph^{011}_1$.

Next, let us consider the two dimensional solutions.  Suppose that
$d$ is as in \refeq{2dderived}, and let
$B=\begin{pmatrix}a_4&a_5\\a_7&a_8\end{pmatrix}$, so $B$ is an
invertible matrix.

Let $d'=\ph^{101}_2+\lambda\ph^{011}_3$. Then setting  $d'g=gd$
leads to the following equalities:
\begin{align*}
lt-nr&=0\\
\lambda(mt-ns)&=0\\
lz-nx&=a_4l+a_5r\\
\lambda(mz-ny)&=a_4m+a_5s\\
0&=a_4n+a_5t\\
rz-tx&=a_7l+a_8r\\
\lambda(sz-ty)&=a_7m+a_8s\\
0&=a_7n+a_8t
\end{align*}
It follows that $B\left[\begin{matrix}n\\t\end{matrix}\right]=0$,
so that $n=0$ and $t=0$.  Then we have $z(ls-mr)=\det(G)\ne 0$,
and
$$
B\begin{pmatrix}l&m\\r&s\end{pmatrix}= z
\begin{pmatrix}l&\lambda m\\r&\lambda s\end{pmatrix}
$$
Note that $z\ne 0$ and $ls-mr\ne 0$. This implies that
$(l,r)$ and $(m,s)$ are linearly independent eigenvectors of $B$
with eigenvalues $z$ and $z\lambda$ respectively.  Moreover, any
matrix $B$ which has linearly independent eigenvectors of this
form, with the ratio of eigenvalues being $\lambda$, determines a
codifferential which is equivalent to $d'$.  In particular, if
$d'_i=\ph^{101}_1+\lambda_i\ph^{011}_2$, for $i=1,2$, then $d'_1$
is equivalent to $d'_2$ when $\lambda_2=\lambda_1$ or
$\lambda_2=1/\lambda_1$.  Thus the nonequivalent solutions of this
form parameterize the punctured unit disc in $\C$, with the
boundary glued together. Moreover, if we consider the case
$\lambda=0$, which is a one dimensional solution, it is then easy
to see that it is equivalent to the solution $d=\ph_1^{110}$.  We
therefore add this solution to the family of two dimensional
solutions.  We shall see later, when
we study versal deformations of these codifferentials,
that this idea makes sense.

Now, if the matrix associated to $d$ is not diagonalizable, then
there is some basis in which its matrix can be expressed in the
form $B=\begin{pmatrix}\lambda&1\\0&\lambda\end{pmatrix}$, so
$d=\ph^{101}_1\lambda+\ph^{101}_2+\ph^{011}_2\lambda$. Let us
consider the special case
$d'=\ph^{101}_1+\ph^{101}_2+\ph^{011}_2$. It is evident that if we
consider the diagonal automorphism $g(\one)=\one$,
$g(\two)=\two\lambda$, $g(\tre)=\tre\lambda$, then $d'g=gd$. Thus
all invertible defective matrices arise from the single
codifferential $d=\ph^{101}_1+\ph^{101}_2+\ph^{011}_2$.

Finally, let us consider the three dimensional case, with matrix
as in \eqref{3ddet}.  If $a_3\ne0$, then by introducing a new
basis element $\tre'=\one a_1+\two a_2 +\tre a_3$, we can replace
$A$ with a simpler matrix of the form
$A=\begin{pmatrix}0&0&1\\a_4& a_5&0\\a_7&-a_4&0\end{pmatrix}$.  By
interchanging the elements $\one$, $\two$ and $\tre$, we can see
that if any one of the elements $a_1$, $a_3$ or $a_7$, does not
vanish, then we can make a similar transformation. On the other
hand, if all of those coefficients vanish, then
$\det(A)=-2a_1a_2a_4$, so none of the three coefficients $a_1$,
$a_2$ and $a_4$ can vanish.  Replacing the basis element $\two$
with $\two'=\two+\tre$ will result in a matrix whose $a_3$ term
does not vanish. Thus we can assume without loss of generality
that $A$ has the simpler form above.

Now if $a_5\ne 0$, then replacing $\two$ with $\one+\two$ results
in a matrix of the form $A=\begin{pmatrix}0&0&1\\0&
1&0\\a_7&0&0\end{pmatrix}$, and a simple linear change of
variables allows one to assume that $a_7=1$. This yields  a
codifferential of the form
$d'=\ph^{110}_3+\ph^{101}_2+\ph^{011}_1$. Otherwise, we have
$a_5=0$, and we may as well assume that $a_7=0$, otherwise, by
interchanging $\one$ and $\two$, we reduce to $d'$.  Now a simple
linear change of variables will allow us to assume $a_4=1$, so our
second candidate for a three dimensional codifferential is
$d=\ph^{110}_3+\ph^{101}_1-\ph^{011}_2$. It is easy to verify that
$G=\begin{pmatrix}-i&0&1\\-\frac
i2&0&-\frac12\\0&i&0\end{pmatrix}$ will yield $d'g=gd$, so these
two codifferentials are equivalent.

Thus, up to equivalence, we obtain precisely the codifferentials
\begin{align*}
d_1&=\ph^{011}_1\\
\dt{\lambda}&=\ph^{101}_1+\ph^{011}_2\lambda\\
d_2&=\ph^{101}_1+\ph^{101}_2+\ph^{011}_2\\
d_3&=\ph^{110}_3+\ph^{101}_2+\ph^{011}_1,
\end{align*}
where $\dt{\lambda}$ is identified with $\dt{\lambda\inv}$.  This
gives one family and three exceptional differentials. However,  as
we shall see when we study versal deformations, this
classification has not yet revealed exactly how the moduli space
fits together.

\section{Versal Deformations of the Codifferentials}
We now begin the process of constructing a miniversal
deformation for each of the codifferentials we have studied.
First, we compute the cohomology of the codifferential, use it to
write the universal infinitesimal deformation, and then apply the
bracket process above to determine a miniversal deformation and
the relations on the parameters.  Along the way, we will discover
that the cohomology of the differentials reveals a lot of
information about the moduli space of three dimensional Lie
algebras.  So far, what we have seen is that there is a family of
codifferentials and three special cases that lie outside the
family.  But how do these special codifferentials fit together
with the family as a moduli space?

\subsection{The Codifferential $d_3=\ph^{110}_3+\ph^{101}_2+\ph^{011}_1$}
This codifferential is the only one with a three dimensional
derived algebra, and corresponds to the simple Lie algebra
$\text{sl}(2,\C)$.  Accordingly, we would expect that the
cohomology vanishes in all dimensions. However, keep in mind that
in the usual definition of Lie algebra cohomology, there is a
space $L_0=\hom(\C,W)=W$ with the bracket of a cochain and an
element of $L_0$ given by $\br{\ph}w(\alpha)=\ph(w\alpha)$.  The
1-cocycles $Z^1$ correspond to derivations of the Lie algebra
structure determined by $d_3$. The 0-coboundaries are the inner
derivations. For a simple Lie algebra, all derivations are inner.

Let us compute a basis of the inner derivations. We have
$$D(\one)(\one)=[\one,d_3](\one)=-[d_3,\one](\one)=d_3(\one\one)=0.$$
Similarly, $D(\one)(\two)=-d_3(\one\two)=-\tre$ and
$D(\one)(\tre)=-\two$. Thus we obtain
\begin{align*}
D(\one)&=-\ph^{010}_3-\ph^{001}_2\\
D(\two)&=\ph^{100}_3-\ph^{001}_1\\
D(\tre)&=\ph^{100}_2+\ph^{010}_1
\end{align*}
From our perspective, since we do not include $L_0$ in our space
of cochains, we expect that $H^1$ will be $3|0$-dimensional, and
all higher cohomology should vanish.  This indeed turns out to be
the case.  We will not reproduce the calculations of the
cohomology here,  but we will give detailed computations later for
some more interesting cases, so the method will be clear.

Nevertheless, since we have three cocycles, all even, we do have a
nontrivial infinitesimal deformation. Let us adopt the convention
to use the Greek letter $\theta$ for an odd parameter, and the
Roman letter $t$ for an even one. For an odd parameter $\theta$,
we have $\theta^2=0$ as a consequence of the graded commutativity,
so we do not consider this to be a relation on our parameter
algebra.  Thus, we have
\begin{equation*}
d_3^1=\ph^{110}_3+\ph^{101}_2+\ph^{011}_1+(\ph^{010}_3+\ph^{001}_2)\theta_1+(\ph^{100}_3-\ph^{001}_1)\theta_2
+(\ph^{100}_2+\ph^{010}_1)\theta_3.
\end{equation*}
In computing $[d_3^1,d_3^1]$, we note that the brackets of the
cocycles with $d_3$ vanish, so we only need to compute the
brackets of the cocycles with each other. Using the fact that the
square of an odd parameter is zero, this means we only need to
calculate the following brackets.
\begin{align*}
\br{\ph^{010}_3+\ph^{001}_2}{\ph^{100}_3-\ph^{001}_1}&=\ph^{100}_2+\ph^{010}_1\\
\br{\ph^{010}_3+\ph^{001}_2}{\ph^{100}_2+\ph^{010}_1}&=\ph^{100}_3-\ph^{001}_1\\
\br{\ph^{100}_3-\ph^{001}_1}{\ph^{100}_2+\ph^{010}_1}&=\ph^{010}_3+\ph^{001}_2
\end{align*}
Thus we obtain
\begin{equation*}
\br{d_3^1}{d_3^1}= (\ph^{100}_2+\ph^{010}_1)\theta_1\theta_2+
(\ph^{100}_3-\ph^{001}_1)\theta_1\theta_3+
(\ph^{010}_3+\ph^{001}_2)\theta_2\theta_3.
\end{equation*}
Of course, these are all cocycles which are not coboundaries, so
we obtain the relations $\theta_i\theta_j=0$ for all $i$, $j$.
Thus the infinitesimal deformation is miniversal, and the base of
the miniversal deformation is given by
$\C[\theta_1,\theta_2,\theta_3]/(\theta_1\theta_2,\theta_1\theta_3,\theta_2\theta_3)$.

Note that the vanishing of $H^2$ is consistent with the
observation that any small change in the codifferential $d_3$ will
give rise to a codifferential $d'$ which will still have a
3-dimensional derived subalgebra.  Thus any small change in $d_3$
gives rise to the same codifferential, and we see that $d_3$ does
not deform into any of the other codifferentials.

\subsection{The Codifferential $d_2=\ph^{101}_1+\ph^{101}_2+\ph^{011}_2$}

First, we give a table of all the coboundaries.
\begin{align*}
D(\ph^{100}_1)&=-\ph^{101}_2,\qquad &D(\ph^{110}_1)&=-\ph^{111}_1+\ph^{111}_2\\
D(\ph^{100}_2)&=0,\qquad &D(\ph^{110}_2)&=-\ph^{111}_2\\
D(\ph^{100}_3)&=\ph^{110}_2+\ph^{101}_3,\qquad &D(\ph^{110}_3)&=-2\ph^{111}_3\\
D(\ph^{010}_1)&=\ph^{101}_1-\ph^{011}_2,\qquad &D(\ph^{101}_1)&=0\\
D(\ph^{010}_2)&=\ph^{101}_2,\qquad &D(\ph^{101}_2)&=0\\
D(\ph^{010}_3)&=-\ph^{110}_1-\ph^{110}_2+\ph^{101}_3+\ph^{011}_3,\qquad &D(\ph^{101}_3)&=\ph^{111}_2\\
D(\ph^{001}_1)&=0,\qquad &D(\ph^{011}_1)&=0\\
D(\ph^{001}_2)&=0,\qquad &D(\ph^{011}_2)&=0\\
D(\ph^{001}_3)&=-\ph^{101}_1-\ph^{101}_2-\ph^{011}_2,\qquad &D(\ph^{011}_3)&=-\ph^{111}_1-\ph^{111}_2\\
\end{align*}
Let $h_n$ be the dimension of $H^n$, the $n$-th cohomology group
of $d$, $b_n$ be the dimension of $B^{n}=D(L_n)$, the space of
$n$-coboundaries, and $z_n$ be the dimension of
$Z^{n}=\ker(D:L_n\ra L_{n+1})$, the space of $n$-cocycles.  Note
that the $n$-coboundaries are a subspace of the $n+1$-cocycles in
this notation.

Let us determine these dimensions for the codifferential $d_2$.
Note that in there are three obvious 1-cocycles, $\ph^{100}_2$,
$\ph^{001}_1$ and $\ph^{001}_2$. In addition we also have the cocycle
$\ph^{100}_1+\ph^{010}_2$. Thus $h_1=z_1=4$, and $b_1=5$.
Evidently, $b_2=3$, so we must have $z_2=6$.  Thus $h_2=1$ and
$h_3=0$. We can choose $\ph^{011}_2$ as a basis of $H^2$. The
universal infinitesimal deformation is given by
\begin{equation*}
d_2^1=\ph^{101}_1+\ph^{101}_2+(1+t)\ph^{011}_2+\ph^{100}_2\theta_1+\ph^{001}_1\theta_2+\ph^{001}_2\theta_4
+(\ph^{100}_1+\ph^{010}_2)\theta_4.
\end{equation*}
The nonzero brackets which occur in the computation of
$[d_2^1,d_2^1]$ are
\begin{align*}
\br{\ph^{011}_2}{\ph^{100}_2}&=\ph^{101}_2,\qquad &\br{\ph^{001}_1}{\ph^{100}_1+\ph^{010}_2}&=-\ph^{001}_1\\
\br{\ph^{100}_1}{\ph^{001}_1}&=\ph^{001}_1,\qquad &\br{\ph^{001}_2}{\ph^{100}_1+\ph^{010}_2}&=-\ph^{001}_2\\
\end{align*}
Thus we have
\begin{equation*}
\br{d_2^1}{d_2^1}=
2\ph^{101}_2t\theta_1+\ph^{001}_2(\theta_1\theta_2-\theta_3\theta_4)-\ph^{001}_1\theta_2\theta_4.
\end{equation*}
Of the three cocycles appearing on the right hand side of this
equation, only the first is a coboundary. Thus we obtain two
second order relations, $\theta_1\theta_2-\theta_3\theta_4=0$ and
$\theta_2\theta_4=0$, and we need to add something to $d_2^1$ to
obtain a second order deformation. Since
$D(\ph^{100}_1)=-\ph^{101}_2$, we can express
$d_2^2=d_2^1+\ph^{100}_1t\theta_1$. We compute
\begin{equation*}
\br{d_2^2}{d_2^2}=
\ph^{001}_2(\theta_1\theta_2-\theta_3\theta_4)-\ph^{001}_1(\theta_2\theta_4+t\theta_1\theta_2).
\end{equation*}
Since no coboundary terms occur, we obtain that $d_2^2$ is a
miniversal deformation, and the base is given by
\begin{equation*}A=
\C[[t,\theta_1,\theta_2,\theta_3,\theta_4]]/(\theta_1\theta_2-\theta_3\theta_4,\theta_2\theta_4+t\theta_1\theta_2)
.\end{equation*} Since there is a nontrivial deformation in the
Lie algebra direction, we can explore how the deformation moves
our codifferential in the moduli space. Note that $d_2^1$ is a
codifferential which has two dimensional derived algebra, and it
has eigenvalues $1$ and $1+t$, so it lies in the family
$\dt{\lambda}$, and is near to $\dt{1}$.  In fact, a punctured
neighborhood of $d_2$ looks exactly like a neighborhood of $\dt1$.
However, $\dt1$ is not close to $d_2$, in the sense that a small
neighborhood of $d_2$ does not contain $\dt1$. When we study
$\dt1$, we shall see that the opposite statement is not true.
\subsection{The Family of Codifferentials $\dt\lambda=\ph^{101}_1+\ph^{011}_2\lambda$}
The cohomology will depend to some extent on the value of
$\lambda$.  The coboundaries are given by

\begin{align*}
D(\ph^{100}_1)&=0,\qquad &D(\ph^{110}_1)&=-\ph^{111}_1\lambda\\
D(\ph^{100}_2)&=\ph^{101}_2(1-\lambda),\qquad &D(\ph^{110}_2)&=-\ph^{111}_2\\
D(\ph^{100}_3)&=\ph^{110}_2\lambda+\ph^{101}_3,\qquad &D(\ph^{110}_3)&=-\ph^{111}_3(\lambda+1)\\
D(\ph^{010}_1)&=-\ph^{011}_1(1-\lambda),\qquad &D(\ph^{101}_1)&=0\\
D(\ph^{010}_2)&=0,\qquad &D(\ph^{101}_2)&=0\\
D(\ph^{010}_3)&=-\ph^{110}_1+\ph^{011}_3\lambda,\qquad &D(\ph^{101}_3)&=\ph^{111}_2\lambda\\
D(\ph^{001}_1)&=0,\qquad &D(\ph^{011}_1)&=0\\
D(\ph^{001}_2)&=0,\qquad &D(\ph^{011}_2)&=0\\
D(\ph^{001}_3)&=-\ph^{101}_1-\ph^{011}_2\lambda,\qquad &D(\ph^{011}_3)&=-\ph^{111}_1\\
\end{align*}
We shall see that the only thing special about the case
$\lambda=0$ is that the dimension of the derived algebra drops to
1, but as far as deformations go, it will behave like
a generic element of the family.  The cases $\lambda=\pm 1$,
however, are not generic in terms of their deformations.  This
makes sense, because in the identification of $td\lambda$ with
$\td{\lambda\inv}$, we see that every point in the unit disc has a
neighborhood that is like a usual disc in $\C$, with the exception
of $\pm 1$, which are orbifold points.  So it is not surprising to
find
 some kind of exceptional behavior for these codifferentials.

\subsubsection{Generic case of $\dt{\lambda}$}
First, we treat the generic case. Clearly, $H^1=\langle
\ph^{100}_1,\ph^{010}_2,\ph^{001}_1,\ph^{001}_2\rangle$, so
$h_1=4$ and $b_1=5$. Evidently, $b_2=3$, so $z_2=6$, $h_2=1$ and
$h_3=0$.  We can choose $\ph^{011}_2$ as a basis of $H^2$. Thus,
the universal infinitesimal deformation is given by
\begin{equation*}
\dt{\lambda}^1=\ph^{101}_1+\ph^{011}_2(\lambda+t)+
\ph^{100}_1\theta_1+\ph^{010}_2\theta_2+\ph^{001}_1\theta_3+\ph^{001}_2\theta_4.
\end{equation*}
It is easy to verify that
\begin{equation*}
\br{\dt{\lambda}^1}{\dt{\lambda}^1}=2\ph^{001}_1\theta_1\theta_3+2\ph^{001}_2\theta_2\theta_4.
\end{equation*}
Thus $\dt{\lambda}^1$ is miniversal and the base of the miniversal
deformation is
$$A=\C[[t,\theta_1,\theta_2,\theta_3,\theta_4]]/(\theta_1\theta_3,\theta_2\theta_4).$$

Looking at the deformation in the Lie algebra direction, we see
that the deformation simply moves along the family.

\subsubsection{The special case $\dt{-1}$}
Now, let us consider the special case $\lambda=-1$. Then we have
an extra cocycle $\ph^{110}_3$ in $H^2$, and correspondingly, an
extra cocycle $\ph^{111}_3$ in $H^3$.  Thus $h_1=4$, $h_2=2$ and
$h_3=1$.  Thus, the universal infinitesimal deformation becomes
\begin{align*}
{\dt{-1}}^1= &\ph^{101}_1+\ph^{011}_2(-1+t_1)+
\ph^{100}_1\theta_1+\ph^{010}_2\theta_2+\ph^{001}_1\theta_3+\ph^{001}_2\theta_4\\
&+\ph^{110}_3t_2+\ph^{111}_3\theta_5.
\end{align*}
Then
\begin{align*}
\frac12\br{\dt{-1}^1}{\dt{-1}^1}&= \ph^{001}_1\theta_1\theta_3+
\ph^{001}_2\theta_2\theta_4+
\ph^{110}_3(t_2\theta_1+t_2\theta_2)\\
&+\ph^{111}_3(\theta_5\theta_1+\theta_5\theta_2-t_1t_2)
-\ph^{111}_1\theta_5\theta_3
-\ph^{111}_2\theta_5\theta_4\\
&+(-\ph^{110}_1-\ph^{011}_3)t_2\theta_3
+(-\ph^{110}_2+\ph^{101}_3)t_2\theta_4.
\end{align*}
Note that the first four terms are cocycles, so they give rise to
second order relations, while the last four terms are
coboundaries, so we need to add corresponding terms to obtain the
second order deformation
\begin{equation*}
\dt{-1}^2=\dt{-1}^1+\ph^{011}_3\theta_3\theta_5+\ph^{110}_2\theta_4\theta_5
-\ph^{010}_3t_2\theta_3-\ph^{100}_3t_2\theta_4.
\end{equation*}
Finally, let us compute the bracket of the second order
deformation with itself. We obtain
\begin{align*}
\frac12\br{\dt{-1}^2}{\dt{-1}^2}&= \ph^{001}_1\theta_1\theta_3+
\ph^{001}_2\theta_2\theta_4+
\ph^{110}_3(t_2\theta_1+t_2\theta_2)\\
&+\ph^{111}_3(\theta_5\theta_1+\theta_5\theta_2-t_1t_2)
+\ph^{011}_3(\theta_3\theta_5\theta_2-t_1t_2\theta_3)\\
&+\ph^{110}_2(\theta_4\theta_5\theta_1-t_1t_2\theta_4)
-\ph^{010}_3t_2\theta_3\theta_2
+\ph^{010}_2t_2\theta_3\theta_4\\
&-\ph^{100}_3t_2\theta_4\theta_1 +\ph^{100}_1t_2\theta_4\theta_3
.\end{align*} All the terms except $\ph^{011}_3$, $\ph^{110}_2$,
$\ph^{010}_3$, and $\ph^{100}_3$ are cocycles, and these
exceptional terms are not even coboundaries.  Thus, by the general
theory, their coefficients must be zero, using the third order
relations
\begin{align*}
\theta_1\theta_3=0,\qquad \theta_2\theta_4=0,\qquad
t_2\theta_3\theta_4=0\\
t_2\theta_1+t_2\theta_2=0,\qquad
\theta_5\theta_1+\theta_5\theta_2-t_1t_2=0
\end{align*}
For example, to see that the coefficient of $\ph^{110}_2$ vanishes, we
observe that
\begin{align*}
\theta_4\theta_5\theta_1-t_1t_2\theta_4&=
\theta_4\theta_5\theta_1-(\theta_5\theta_1+\theta_5\theta_2)\theta_4\\
&=\theta_4\theta_5\theta_1-\theta_5\theta_1\theta_4=0.
\end{align*}

Let us consider the induced topology on the moduli space of
equivalence classes of codifferentials. This topology is not
Hausdorff.  If every neighborhood of a point $a$ contains the
point $b$, then we note that $a$ is in the closure of $b$.  In
this case, we shall say that $a$ is {\it infinitesimally close} to $b$.
Note that $\dt{-1}$ is infinitesimally close to $d_3$, but not the
other way around.  In some sense, this explains the extra
infinitesimal deformation in the Lie algebra direction.

Now at first it may seem strange that this is the first case where
one of our codifferentials deforms into $d_3$.  After all,
generically, we would expect the matrix of a codifferential to be
invertible.  But look carefully at \refeq{3dderived}, and the form
of a solution in the family, and it becomes clear that only when
$\lambda=-1$ can a small change in the codifferential give a
solution satisfying \refeq{3dderived}.

Notice also, that if we neglect odd parameters, the versal
deformation has relation $t_1t_2=0$.  In the classical sense, this
means that if we were to consider an infinitesimal deformation of
the form $d'=d+t(a\ph^{011}_2+b\ph^{110}_3$, then the only cases
where this deformation extends to a second order deformation is
when either $a=0$ or $b=0$, in which case, the infinitesimal
deformation extends trivially.
\subsubsection{The special case $\dt1$}
In this case, there are two additional 1-cocycles, $\phi^{100}_2$
and $\phi^{010}_1$.  Thus $h_1=6$, and $b_1=3$.  Since $b_2=3$, we
see that $h_2=3$ and $h_3=0$.  So we pick up two extra 2-cocycles,
$\ph^{101}_2$ and $\ph^{011}_1$.  It is convenient to replace the
cocycle $\ph^{011}_2$, which we used as a basis of the cohomology
in the generic case, with $\ph^{101}_1-\ph^{011}_2$, because it
simplifies the interpretation of the bracket of the universal
infinitesimal deformation with itself. Thus,
\begin{align*}
\dt{1}^1&=\ph^{101}_1+\ph^{011}_2\lambda
+(\ph^{101}_1-\ph^{011}_2)t_1
+\ph^{101}_2t_2+\ph^{011}_1t_3\\
&+\ph^{100}_1\theta_1+\ph^{010}_2\theta_2+\ph^{001}_1\theta_3+\ph^{001}_2\theta_4
+\ph^{100}_2\theta_5+\ph^{010}_1\theta_6,
\end{align*}
and we compute that
\begin{align*}
\frac12\br{\dt{1}^1}{\dt{1}^1}&=
(\ph^{101}_1-\ph^{011}_2)(t_3\theta_5-t_2\theta_6)
-\ph^{101}_2(2t_1\theta_5+t_2(\theta_2-\theta_1))\\
&+\ph^{011}_1(2t_1\theta_6+t_3(\theta_2-\theta_1))
-\ph^{100}_1\theta_5\theta_6
+\ph^{010}_2\theta_5\theta_6\\
&+\ph^{001}_1(\theta_1\theta_3-\theta_4\theta_6)
+\ph^{001}_2(\theta_2\theta_4-\theta_3\theta_5)\\
&+\ph^{100}_2(\theta_2\theta_5-\theta_1\theta_5)
+\ph^{010}_1(\theta_1\theta_6-\theta_2\theta_6),
\end{align*}
which is precisely the set of cocycles appearing in the universal
infinitesimal deformation, multiplied by the second order
relations.  Thus the infinitesimal deformation is miniversal.

Now let us interpret how $\dt1$ fits into the moduli space.  Note
that there is a deformation along the family, given by the cocycle
$\ph^{101}_1-\ph^{011}_2$, and two other directions of
deformation, each of which corresponds to a deformation of $\dt1$
into the special codifferential $d_2$.  In fact, if we consider
the three dimensional deformation space parameterized by
$(t_1,t_2,t_3)$,  we see that precisely two curves correspond to a
deformation in the $d_2$ direction.  Thus we see that $\dt1$ is
infinitesimally close to $d_2$, although the converse is not true.

So far, we have been able to associate two of the three special
codifferentials with the family in some manner.
\subsection{The Codifferential $d_1=\ph^{011}_1$}
Now we come to the most complicated of the codifferentials. This
is not surprising, because $d_1$ gives a nilpotent Lie algebra
structure, so that we expect to find a lot of deformations. The
table of codifferentials is given by
\begin{align*}
D(\ph^{100}_1)&=\ph^{011}_1,\qquad &D(\ph^{110}_1)&=0\\
D(\ph^{100}_2)&=-\ph^{101}_1+\ph^{011}_2,\qquad &D(\ph^{110}_2)&=\ph^{111}_1\\
D(\ph^{100}_3)&=\ph^{110}_1+\ph^{011}_3,\qquad &D(\ph^{110}_3)&=0\\
D(\ph^{010}_1)&=0,\qquad &D(\ph^{101}_1)&=0\\
D(\ph^{010}_2)&=-\ph^{011}_1,\qquad &D(\ph^{101}_2)&=0\\
D(\ph^{010}_3)&=0,\qquad &D(\ph^{101}_3)&=\ph^{111}_1\\
D(\ph^{001}_1)&=0,\qquad &D(\ph^{011}_1)&=0\\
D(\ph^{001}_2)&=0,\qquad &D(\ph^{011}_2)&=0\\
D(\ph^{001}_3)&=-\ph^{011}_1,\qquad &D(\ph^{011}_3)&=0\\
\end{align*}
Clearly, $h_1=6$, so $b_1=3$.  Since $b_2=1$, we have $h_2=5$ and
$h_3=2$. Before constructing a versal deformation, let us
analyze how this codifferential sits in the moduli space.  Clearly
there are a lot of directions one can deform. Notice that because
$\ph^{101}_1-\ph^{011}_2$ and $\ph^{110}_1+\ph^{011}_3$ are
coboundaries, there remain three ways to deform $d_1$ into $d_3$,
via the cocycles $\ph^{110}_2$, $\ph^{101}_3+\ph^{101}_2$ and
$\ph^{101}_3-\ph^{101}_2$. Generically, if we add small multiples
of these cocycles, we will obtain a codifferential which is
equivalent to $d_3$. Thus $d_1$ is infinitesimally close to $d_3$.

Next, if we add an appropriate multiple of
$\ph^{101}+a\ph^{101}_2$, then we are constructing an element in
the family, with a matrix $B$ given by
$B=\begin{pmatrix}t&at\\1&0\end{pmatrix}$, whose eigenvalues are
given by $\lambda_{\pm}=\frac{t+\sqrt{t^2+ta}}2$.  Set $a=ct$.
Then
\begin{equation*}
\lambda=\frac{\lambda_-}{\lambda_+}=\frac{1-\sqrt{1+4c}}{1+\sqrt{1+4c}}
\end{equation*}
determines the element of the family.  Solving for $c$ in terms of
$\lambda$ we obtain $c=\frac{-\lambda}{(\lambda+1)^2}$, which
gives an element of the family for any value of $\lambda$ except
$\lambda=\pm 1$. The reason that we do not obtain an element of
the family for $\lambda=1$ is that the resulting matrix is
defective, so we obtain $d_2$ instead. There is a solution for
$\lambda=-1$, however, given by the cocycle $\ph^{101}_2$. Thus we
see that $d_1$ is infinitesimally close to $d_2$ and to every member
of the family except for $\lambda=1$.  Moreover, it is not hard to
check that if we take $d$ to be of the form
\begin{equation*}
d=\ph^{011}_1+\ph^{110}_1t_1+\ph^{101}_1t_2+(\ph^{110}_2-\ph^{101}_3)t_3+\ph^{110}_3t_4+\ph^{101}_2t_5,
\end{equation*}
then there is no automorphism which takes it to $\dt1$.

Thus we conclude that $d_1$ is infinitesimally close to every
codifferential except $\dt1$. From the behavior of the elements
$d_2$ and $\dt1$, it is more natural to consider $d_2$ as a member
of the family, because it behaves more consistently with the other
members of the family than $\dt1$. Note, for example, that
$h_2(d_2)=1$, the same as the generic elements of the family,
while $h_2(\dt1)=3$.

The universal infinitesimal deformation is given by
\begin{align*}
\d_1^1&=
\ph^{011}_1+\ph^{110}_1t_1+\ph^{101}_1t_2+(\ph^{110}_2-\ph^{101}_3)t_3+\ph^{110}_3t_4+\ph^{101}_2t_5\\
&+(\ph^{100}_1+\ph^{010}_2)\theta_1+(\ph^{010}_2-\ph^{001}_3)\theta_2
+\ph^{010}_1\theta_3
+\ph^{010}_3\theta_4\\
&+\ph^{001}_1\theta_5 +\ph^{001}_2\theta_6 +\ph^{111}_2\theta_7
+\ph^{111}_3\theta_8
\end{align*}
Unfortunately, the bracket of $d_1^1$ with itself has terms
involving almost every cocycle and coboundary.
\begin{align*}
\br{d_1^1}{d_1^1}&=
+\ph^{110}_1(t_1\theta_1+t_1\theta_2+t_2\theta_4)
+\ph^{101}_1(t_1\theta_6-t_2\theta_2)\\
&+(\ph^{110}_2-\ph^{101}_3)(t_3\theta_1-t_4\theta_6+t_5\theta_4)
+2\ph^{110}_3(-t_3\theta_4+t_4\theta_1+t_4\theta_2)\\
&+2\ph^{101}_2(t_3\theta_6-t_5\theta_2)
-(\ph^{010}_2-\ph^{001}_3)\theta_4\theta_6
-\ph^{010}_1(\theta_2\theta_3+\theta_4\theta_5)
\\&-\ph^{010}_3(\theta_1\theta_4+2\theta_2\theta_4)
+\ph^{001}_1(\theta_1\theta_5+\theta_2\theta_5+\theta_3\theta_6)
\\&+\ph^{001}_2(\theta_1\theta_6+2\theta_2\theta_6)
+\ph^{111}_2(t_1t_5-t_2t_3-\theta_1\theta_7+\theta_2\theta_7+\theta_6\theta_8)
\\&+\ph^{111}_3(-t_1t_3-t_2t_4-2\theta_1\theta_8-\theta_2\theta_8+\theta_4\theta_7)
\\&+(-\ph^{101}_1+\ph^{011}_2)(-t_3\theta_5+t_5\theta_3)
-(\ph^{110}_1+\ph^{011}_3)(t_3\theta_3+t_4\theta_5)
\\&+\ph^{011}_1(-t_1\theta_5+t_2\theta_3)+\ph^{111}_1(\theta_3\theta_7+\theta_5\theta_8)
\end{align*}

If we forget about all the $\theta$ terms, then the second order
relations reduce to $t_1t_5-t_2t_3=0$ and $-t_1t_3-t_tt_5=0$.  To
see why we should expect such relations,  consider the  matrix
$A=\begin{pmatrix}t_1&t_3&t_4\\t_2&t_5&-t_3\\1&0&0\end{pmatrix}$
associated to the 2-cochain part of $d_1^1$. We know that if its
determinant is nonzero, then it must be of the form
\eqref{3ddet}. But then $t_1=t_2=0$.  Otherwise, as the relations
yield, we must have dependent columns in the matrix.

Note that the last four terms are coboundaries, so we need to add
something to $d_1^1$ to obtain a second order deformation.  We
may choose
\begin{align*}
d_1^2&=d_1^2 +\ph^{100}_2(-t_3\theta_5+t_5\theta_3)
+\ph^{100}_3(-t_3\theta_3-t_4\theta_5)
+\ph^{100}_1(-t_1\theta_5+t_2\theta_3)
\\&-\ph^{110}_2(\theta_3\theta_7+\theta_5\theta_8).
\end{align*}
Let
\begin{align*}
b_1&=-\theta_3\theta_7-\theta_5\theta_8\\
b_2&=-t_3\theta_5+t_5\theta_3\\
b_3&=-t_3\theta_3-t_4\theta_5\\
b_4&=-t_1\theta_5+t_2\theta_3
\end{align*}
Then
\begin{align*}
\frac12\br{d_1^2}{d_1^2}&=
+\ph^{110}_1(t_1\theta_1+t_1\theta_2+t_2\theta_4-b_1\theta_3)
+\ph^{101}_1(t_1\theta_6-t_2\theta_2)\\
&+(\ph^{110}_2-\ph^{101}_3)(t_3\theta_1-t_4\theta_6+t_5\theta_4)
\\&+\ph^{110}_3(-2t_3\theta_4+2t_4\theta_1+2t_4\theta_2+t_4b_4-t_1b_3-b_1\theta_4)\\
&+\ph^{101}_2(2t_3\theta_6-2t_5\theta_2+t_5b_4-t_2b_2+b_1\theta_6)
-(\ph^{010}_2-\ph^{001}_3)\theta_4\theta_6
\\&-\ph^{010}_1(\theta_2\theta_3+\theta_4\theta_5-b_4\theta_3)
\\&-\ph^{010}_3(\theta_1\theta_4+2\theta_2\theta_4)
+\ph^{001}_1(\theta_1\theta_5+\theta_2\theta_5+\theta_3\theta_6+b_4\theta_5)
\\&+\ph^{001}_2(\theta_1\theta_6+2\theta_2\theta_6+b_2\theta_5)
\\&+\ph^{111}_2(t_1t_5-t_2t_3-\theta_1\theta_7+\theta_2\theta_7+\theta_6\theta_8-b_4\theta_7-t_2b_1)
\\&+\ph^{111}_3(-t_1t_3-t_2t_4-2\theta_1\theta_8-\theta_2\theta_8+\theta_4\theta_7-b_4\theta_8)
\\&+\ph^{100}_1(-b_3\theta_5-b_2\theta_3)
+\ph^{010}_2b_2\theta_3 +\ph^{001}_3b_3\theta_5
\\&+\ph^{110}_2(t_3b_4-t_1b_2+b_1b_4+b_1\theta_1)
+\ph^{100}_2(-b_3\theta_6+b_2b_4-b_2\theta_2)
\\&+\ph^{101}_3(-t_3b_4-t_2b_3)
-\ph^{011}_2b_1\theta_5
\\&+\ph^{100}_3(b_3b_4+b_3\theta_2+b_3\theta_1-b_2\theta_4)
\end{align*}
In the end, after some work, one sees that the only remaining
coboundary term is
$(-\ph^{101}_1+\ph^{011}_2)(\theta_3\theta_5\theta_7)$, so we get
a third order deformation
\begin{equation*}d_1^3=d_1^2+\ph^{100}_2\theta_3\theta_5\theta_7.\end{equation*}
A bit more work, and we can conclude that this is as far as we
need to go.  We won't give the explicit relations on the base
here.  The main point is that the process stops at $d_1^3$, and
this gives a miniversal deformation of $d_1$.  Note that in our
calculations, using the general theory,  it is safe to ignore
terms which arise in the bracket $d^i$ with itself that do not
contribute to cocycles or coboundaries, because they have a
coefficient which must vanish. Nevertheless, it is possible to
check that the coefficients of the non-cocycle terms in the
bracket of $d_2^2$ with itself, do vanish, using the third order
relations, with the exception of some fourth order terms that
cancel after adding the new fourth order terms arising in the
bracket of $d_2^3$ with itself. Nevertheless, checking the
coefficients is a useful method of avoiding misprints in the
terms.

\subsection{Deformations of the Trivial Codifferential $d_0=0$}
There is one case which we have not touched on, the case of the
trivial codifferential $d=0$.  At first it may seem as if there is
little to obtain from looking at this situation, because it
evidently must deform into every possible type of codifferential,
so we know that it is infinitesimally close to every point in the
moduli space.  Moreover, since there are no coboundaries, the
infinitesimal deformation is obviously miniversal.

On the other hand, we do obtain some relations, and these
relations tell us something about how the moduli space is put
together.  In addition, all second order relations can be
determined from the relations on the zero codifferential, by using
appropriate values for the coefficients. Keep in mind that no
information on higher order relations can be obtained in this
manner.

We do not need a table of coboundaries for $d=0$. Every cochain is
a cocycle, so the universal infinitesimal deformation is given by
\begin{align*}d_0^1&=
\ph^{110}_1t_1+\ph^{110}_2t_2+\ph^{110}_3t_3
+\ph^{101}_1t_4+\ph^{101}_2t_5+\ph^{101}_3t_6
\\&+\ph^{011}_1t_7+\ph^{011}_2t_8+\ph^{011}_3t_9
+\ph^{100}_1\theta_1+\ph^{100}_2\theta_2+\ph^{100}_3\theta_3
\\&+\ph^{010}_1\theta_4+\ph^{010}_2\theta_5+\ph^{100}_3\theta_6
+\ph^{001}_1\theta_7+\ph^{001}_2\theta_8+\ph^{001}_3\theta_9
\\&+\ph^{111}_1\theta_{10}+\ph^{111}_2\theta_{11}+\ph^{111}_3\theta_{12}
.
\end{align*}
We will not give a table of all the relations that are obtained
here, but from such a table, one can construct the second order
relations for any codifferential by simply substituting the
$\theta$'s and $t$'s in the coefficient of a term in the bracket
with the coefficients that occur in that particular infinitesimal
deformation in the same places.
\section{Conclusion}
The main purpose of this paper was to illustrate how to compute
versal deformations of Lie algebras.  In fact, we worked in a
slightly more general picture, that of \linf\ algebras.  The
classification of Lie algebras of dimension three is well known,
but we think that by studying the deformations more closely,
the picture of the geometry of the moduli space of Lie algebra
structures on a three dimensional vector space becomes much
clearer.

\providecommand{\bysame}{\leavevmode\hbox to3em{\hrulefill}\thinspace}
\providecommand{\MR}{\relax\ifhmode\unskip\space\fi MR }
\providecommand{\MRhref}[2]{%
  \href{http://www.ams.org/mathscinet-getitem?mr=#1}{#2}
}
\providecommand{\href}[2]{#2}


\begin{thebibliography}{10}

\bibitem{aga}
Y.~Agaoka, \emph{On the variety of 3-dimensional {L}ie algebras}, Lobachevskii
  Journal of Mathematics \textbf{3} (1999), 5--17.

\bibitem{bfls}
G.~Barnich, R.~Fulp, T.~Lada, and J.~Stasheff, \emph{The sh {L}ie structure of
  {P}oisson brackets in field theory}, Comm. Math. Phys. \textbf{191} (1998),
  no.~3, 585--601.

\bibitem{bfp1}
D.~Bodin, A.~Fialowski and M.~Penkava,
{\it Classification of $L_\infty$ Algebras on a $2|1$-dimensional Space}
(preprint 2003).

\bibitem{fi}
A.~Fialowski, \emph{Deformations of {L}ie algebras}, Mathematics of the
  USSR-Sbornik \textbf{55} (1986), no.~2, 467--473.

\bibitem{ff3}
A.~Fialowski and D.~Fuchs, \emph{Singular deformations of {L}ie algebras on an
  example}, Topics in Singularity Theory (Providence, RI) (A.~Varchenko and
  V.~Vassilie, eds.), A.M.S. Translation Series 2, Vol.180, Amer.\ Math.\ Soc.,
  1997, V.~I.~Arnold {\rm 60}${}^{th}$ Anniversary Collection,.

\bibitem{ff2}
\bysame, \emph{Construction of miniversal deformations of {L}ie algebras},
  Journal of Functional Analysis (1999), no.~161(1), 76--110.

\bibitem{fp1}
A.~Fialowski and M.~Penkava, \emph{Deformation theory of infinity algebras},
  Journal of Algebra \textbf{255} (2002), no.~1, 59--88, math.RT/0101097.

\bibitem{fp2}
\bysame, \emph{Examples of infinity and {L}ie algebras and their versal
  deformations}, Geometry and Analysis on {L}ie groups, Banach Center
  Publications, {\bf 55} (2002), pp.~27--42, math.QA/0102140.

\bibitem{fp4}
\bysame, \emph{Classification and Extensions of \linf\ Algebras of
  Dimension $1|2$}, preprint, 2003.

\bibitem{jac}
Nathan Jacobson, \emph{Lie algebras}, John Wiley \& Sons, 1962.

\bibitem{lm}
T.~Lada and M.~Markl, \emph{Strongly homotopy {L}ie algebras}, Comm. in Algebra
  \textbf{23} (1995), 2147--2161.

\bibitem{ls}
T.~Lada and J.~Stasheff, \emph{Introduction to sh {L}ie algebras for
  physicists}, Preprint hep-th 9209099, 1990.

\bibitem{ss}
M.~Schlessinger and J.~Stasheff, \emph{The {L}ie algebra structure of tangent
  cohomology and deformation theory}, Journal of Pure and Applied Algebra
  \textbf{38} (1985), 313--322.

\bibitem{sta3}
J.D. Stasheff, \emph{Closed string field theory, strong homotopy {L}ie algebras
  and the operad actions of moduli spaces}, Perspectives in Mathematical
  Physics, Internat. Press, Cambridge, MA, 1994, pp.~265--288.

\bibitem{tu}
H.~Tasaki and M.~Umehara, \emph{An invariant on 3-dimensional {L}ie algebras},
  Proceedings of the American Mathematical Society \textbf{115} (1992), no.~2,
  293--294.

\end{thebibliography}
\end{document}